\def\IR{{\Bbb R}} 
 
\def\IS{{\Bbb S}} 
\def\IK{{\Bbb K}}

\def\IC{\Bbb C} 
\def\ID{{\Bbb D}}

\def\tOmega{\tilde{\Omega}}
\def\zbar{{\overline{z}}} 
\def\zetabar{{\overline{\zeta}}} 
\def\wbar{{\overline{w}}}

\documentclass[12pt]{article} 

\usepackage[psamsfonts]{amssymb} 
\usepackage{amsfonts,amsmath}

\usepackage{epsfig,multicol}

\newtheorem{theorem}{Theorem}
\newtheorem{lemma}{Lemma}
\newtheorem{corollary}{Corollary}

\title{Harmonic degree $1$ maps are diffeomorphisms: Lewy's theorem for curved metrics.} 
\author{Gaven J. Martin \thanks{Research supported in part by 
grants from the  N.Z. Marsden Fund. \newline  
NZ Inst. for Advanced Study,  Massey University \& Magdalen College, Oxford University  \newline
G.J.Martin@massey.ac.nz}
} 
\date{} 
\begin{document}

\maketitle 


\begin{abstract}  In 1936 H. Lewy showed that the Jacobian determinant of a harmonic homeomorphism between planar domains does not vanish and thus the map is a diffeomorphism.  This built on the earlier existence results of Rad\'o and Kneser.  R. Shoen and S.T. Yau generalised this result to degree $1$ harmonic mappings between closed Riemann surfaces.  Here we give a new approach that establishes all these results  in complete generality.
\end{abstract}

\section{Introduction} 

In 1936 H. Lewy \cite{Lewy} showed that the Jacobian determinant of a harmonic homeomorphism between planar domains does not vanish.  This built on the earlier existence results of T. Rad\'o and H. Kneser, \cite{Rado,Kneser} and subsequently G. Choquet \cite{Choquet}.  A modern framework would place Lewy's result as establishing an optimal degree of regularity for minimisers of a variational mapping problem - DIrichlet energy minimising homeomorphisms - and in this framework there is a considerable literature and  there have been attempts to generalise Lewy's result in different directions.  See for instance \cite{AN,AIMO, Pavlovich, AIM} and the references therein to start with.

Earlier,  in 1978,  R. Schoen and S.T. Yau proved that a harmonic homeomorphism (actually they need only degree $1$) between closed Riemann surfaces of negative curvature is a diffeomorphism, \cite{SY}.   Their proof was global in nature and they needed boundary conditions in order to relax the assumption that the surface was closed.

There is a substantial literature on the relationship between quasiconformal and harmonic mappings.  This is because of the relationships between the various models of Teichm\"uller space and the metrics on these moduli spaces,  see M. Wolf \cite{Wolf} and T. Wan \cite{Wan} to start with. Most of this literature considers only quasiconformal harmonic diffeomorphisms.   

Here we will establish Lewy's theorem for degree $1$ mappings between planar domains which are harmonic with respect to an arbitrary metric,  Theorem \ref{MR}. When the metric is of either non-negative curvature or of non-positive curvature,  we obtain interesting ancillary results concerning the maximum/minimum principles for the distortion of these mappings.  Our approach will further refine the connections between quasiconformal and harmonic mappings.

\medskip

The definition of a harmonic mapping does not depend on the conformal structure of the domain of definition.  We   therefore present our main theorem in the following way.  Precise definitions are given below.

\begin{theorem} \label{MR} Let $\Omega, \tOmega\subset \IC$ be planar domains and $\rho\in C^\infty(\tOmega,\IR_+)$  defining a metric $\rho(z)|dz|$.  Suppose $f:\Omega\to(\tOmega,\rho)$ is a harmonic mapping with respect to the metric $\rho$.  Let $U\subset\Omega$ be relatively compact.

Then there is a quasiconformal diffeomorphism $g: U  \to \IC$ and a holomorphic function $\varphi:g(U) \to \IC$ such that $f = \varphi \circ g$.  In particular if $f$ has degree $1$,  then  $f$ is a diffeomorphism.

Further,  let ${\cal K}_\rho$ denote the curvature of the metric $\rho(z)|dz|$. Then
\begin{itemize}
\item ${\cal K}_\rho\geq 0$ implies the distortion of $g$ is subhamonic and cannot achieve a maximum on any relatively compact set,  and
\item  ${\cal K}_\rho\leq 0$ implies the distortion of $g$ is superhamonic and cannot achieve a minimum on any relatively compact set,  unless that minimum is $0$.
\end{itemize}
\end{theorem}
 
In general, the restriction to a subdomain $U$ is necessary to guarantee quasiconformality even for harmonic self homeomorphism of the unit disk, \cite{Pavlovich}.  

Our proof of this result proceeds in the following manner.  First,  we show that a locally quasiconformal harmonic mapping is a diffeomorphism.  This is basically via a degree argument.  Then we show that the Beltrami coefficient $\mu_f$ of a harmonic mapping has $|\mu_f|\in C^2(\Omega)$.  This is based around an analysis of the power series expansion of a harmonic mapping.  Finally we analyse  $\log |\mu_f|$ establishing a ``super-averaging'' property  (and  sub/super hamonicity in various cases) via a well know differential inequality of Heinz.  This shows the mapping $f$ is locally quasiregular.  Finally  the result will follow from the Stoilow factorisation theorem,  giving the holomorphic $\varphi$,  and the first part.

\subsection{Harmonic mappings}
Let $\Omega, \tOmega\subset \IC$ be planar domains and $\rho:\tOmega\to\IR_+$ a non-vanishing $C^\infty(\tOmega)$ function.    
 A   mapping $f:\Omega \to(\tOmega,\rho)$ is said to be harmonic with respect to the metric $\rho(z)|dz|$ if it is a solution to  the nonlinear second order elliptic equation, called the {\em tension equation}, 
\begin{equation}\label{te}
f_{z\zbar}(z) +   (\log \rho)_z(f)\; f_z(z)\; f_\zbar(z) = 0.
\end{equation} 
This equation is the Euler-Lagrange equation for the variational integral,  defined for Lipschitz mappings,  
\begin{equation}
{\cal E}(f)=\int_\Omega \|\nabla f\|^2\, \rho^2(f) \, dz
\end{equation}
and which is   the Dirichlet integral should $\rho \equiv 1$.  We remark that given our assumptions on $\rho$,  a solution $f$ to (\ref{te}) is a smooth mapping,  $f\in C^\infty(\Omega,\tOmega)$.  Further,  because the coefficient $(\log \rho)_z(f)$  of (\ref{te}) depends on $f$ in an extremely nonlinear way,  much of the theory of elliptic equations cannot be utilised and existence is a very interesting question which is unsettled for the most interesting cases of mappings between the hyperbolic disk with prescribed boundary values - the Schoen conjecture being a focus here. 
 
\subsection{Quasiconformality}

In what follows,  the Beltrami coefficient of a quasiconformal mapping will play a crucial role.  Given a mapping $f$ between planar domains with Sobolev regularity $f\in W^{1,2}_{loc}(\Omega,\IC)$ we define the  the complex dilatation (or Beltrami coefficient) of $f$ as the measurable complex valued  function
\begin{equation}
\mu_f(z) = \frac{f_\zbar(z)}{f_z(z)}
\end{equation}
We will always assume that $f$ preserves orientation so that the Jacobian determinant is non-negative,  $J(z,f)=|f_z|^2-|f_\zbar|^2\geq 0$,  and  hence 
\begin{equation}
|\mu_f(z)|\leq 1
\end{equation}
A mapping $f\in W^{1,2}_{loc}(\Omega,\IC)$   is said to be {\em quasiregular} if 
\begin{equation}
\|\mu_f(z)\|_{L^\infty(\Omega)} = k < 1
\end{equation}
and {\em quasiconformal} if,  in addition,  the mapping is a homeomorphism.
The distortion of a quasiconformal or quasiregular mapping is the function
\begin{equation} 
\IK(z,f) = \frac{1+|\mu_f(z)|^2}{1-|\mu_f(z)|^2},  \hskip20pt \IK =\|\IK\|_{L^\infty(\Omega)}  = \frac{1+k^{2}}{1-k^{2}}
\end{equation}
It follows that $f$ satisfies the distortion inequality between the differential matrix $\nabla f$ and its Jacobian determinant.
\[ \|\nabla f\|^2 \leq \IK \, J(z,f) \]
The modern theory of planar quasiconformal mappings is accounted in \cite{AIM}.

\subsection{The Hopf differential}

We define the Hopf differential of a harmonic mapping $f$ to be
\begin{equation}
\Phi_f(z) = \rho^2(f(z)) f_z(z) \overline{f_\zbar(z)}.
\end{equation}
The following elementary computation shows that the well known fact that the Hopf differential is holomorphic.
 \begin{eqnarray*}
\Phi_\zbar(z) & = &  [\rho^2(f)]_\zbar f_z \overline{f_\zbar} +  \rho(f) f_{z\zbar}  \overline{f_\zbar} +  \rho(f) f_z \overline{f_{z\zbar}} \\
& = &  2\rho(f)\,\Big( [\rho(f)]_\zbar f_z \overline{f_\zbar}   -   \rho_z(f)  f_zf_\zbar \overline{f_\zbar}  -     \overline{ \rho_z(f) f_zf_\zbar} f_z\Big) \\
& = & 2\rho(f)f_z \overline{f_\zbar}  \,\Big(  [\rho(f)]_\zbar  -   \rho_z(f)    f_\zbar   -     \rho_\zbar(f)  \overline{f_z} \Big) = 0
\end{eqnarray*}
Indeed one can make this calculation at a distributional level to see $\Phi_f$ is holomorphic as soon as it lies in $W^{1,1}_{loc}(\ID)$ by the Looman-Menchoff theorem.  However there are rather stronger refinements,  see \cite{Iwaniecetc}.   

 \medskip
 
If the holomorphic Hopf differential $\Phi_f$ is not identically zero,  then it has isolated zeros of finite order.  We will use this fact repeatedly.
 
 \begin{corollary} A harmonic function $f:\Omega\to(\Omega,\rho)$ has zeros of finite order and both the sets  $\{z:f_z(z)=0\}$ and $\{z:f_\zbar(z)=0\}$ are discrete in  $\Omega$.  Further,  since $f$ is orientation preserving
\[\{z:f_\zbar(z)=0\} \subset  \{z:f_z(z)=0\}.\]
 \end{corollary}

 \section{Harmonic maps and complex dilatation}
 
We now establish that quasiconformal harmonic mappings are diffeomorphisms.
 
 \begin{theorem} \label{hqc} Let $\rho\in C^{\infty}(\Omega)$ be a smooth positive density and $f:\Omega\to (\tOmega,\rho)$ a harmonic quasiconformal mapping.  Then $\nabla f \neq 0$ and
 $\mu_f \in C^{\infty}(\ID)$.
 \end{theorem}
 This result follows directly from the rather more precise result below.
  \begin{theorem}  Let $\rho\in C^{k+1}(\Omega)$ be a positive density and $f:\Omega\to (\tOmega,\rho)$ a harmonic $\IK$-quasiconformal mapping.  Suppose that
\begin{equation}\label{reghyp}
k >   \IK-1+\sqrt{\IK^2-1}
\end{equation} 
  Then $J(z,f) > 0 $ on $\Omega$ and $ \mu_f \in C^{k}(\Omega)$.
 \end{theorem}
 \noindent{\bf Proof.} Since
\[ J(z,f) = |f_z|^2 - |f_\zbar|^2 = (1-|\mu_f|^2) |f_z|^2  \]
we need only show $f_z\neq 0$.  Note that it is an elementary consequence of the regularity of the inhomogeneous Laplace equation $\Delta f = u$ and the structure of the tension equation (\ref{te}), that $f$ basically has two more derivatives than $(\log \rho)_z \in C^{k}(\ID)$.  Our arguments can therefore be further refined as we shall only  assume $f\in C^{k+1}(\ID)$.  By pre and post-composition with  translations we may assume that $f(0)=0$ and seek to show $f_z(0)\neq 0$.  Since $f$ is $\IK$-quasiconformal, so is its inverse and therefore $f$  is $\alpha =( \IK-\sqrt{\IK^2-1})\leq 1$ bi-H\"older continuous,  see \cite[Theorem 3.10.2]{AIM}.  In particular we have the estimate
 \[ c_0 |z|^{1/\alpha} \leq |f(z)| \leq C_0 |z|^\alpha \]
 near the origin.  We now make an argument based on the Taylor series for $f$ which we will significantly refine later.  The bound on the left-hand side here, together with  Taylor's formula,   implies that for some smallest integer $M \leq 1/\alpha = \IK+\sqrt{\IK^2-1}$ and $p+q=M$,
\begin{equation}\label{8}  \Big(\frac{\partial^M\, f}{\partial z^p \partial \zbar^q}\Big)(0) \neq  0
\end{equation}
since the regularity hypothesis (\ref{reghyp}) assures us that $f\in C^{M}(\ID)$.   Actually here, in the $C^\infty$ case, we could argue this from the fact that $0$ can be a zero of only finite order. 
The choice of $M$ as smallest implies 
\begin{equation}\label{9}  \Big(\frac{\partial^{s}\, f}{\partial z^{s} }\Big)(0) =  \Big(\frac{\partial^{s}\, f}{\partial \zbar^{s} }\Big)(0) = 0, \hskip10pt s<M.
\end{equation}
However,  repeatedly differentiating the tension equation shows that if $p\neq M$ or $q\neq M$
\[ \Big(\frac{\partial^M\, f}{\partial z^p \partial \zbar^q}\Big)(0) =  0 \]
Thus 
\[ a =  \Big(\frac{\partial^M\, f}{\partial z^M  }\Big)(0) \neq   0, \hskip10pt {\rm or } \hskip10pt  b= \Big(\frac{\partial^M\, f}{ \partial \zbar^M}\Big)(0) \neq  0 \]
and $J(0,f)\geq 0$ implies $|a|\geq |b|$. We deduce that there is some $1 < m \leq \IK-1+\sqrt{\IK^2-1}$ such that  $f$ has the Taylor series expansion near $0$,
 \begin{equation}
 f(z) = a z^m + b \zbar^m  + o(|z|^m), \hskip20pt |a|\geq |b|.
 \end{equation}
 We next compute
\begin{eqnarray*} L_r= \max_{|z|=r} |f(z)-f(0)| &=& (|a|+|b|)r^m+o(r^m) \\
 \ell_r = \min_{|z|=r}   |f(z)-f(0)|& =& \big| |a|-|b| \big| r^m+o(r^m)
 \end{eqnarray*}
  Now, as $f$ is quasiconformal,   Mori's distortion theorem \cite[\S 3.10]{AIM} establishes an upper bound on $L_r/\ell_r$ which is independent of  $r$ (and depends only on $\IK$). Thus
  \begin{equation}\label{13}
  |a| >  |b|
  \end{equation} 
Next,  for $r>0$ but small enough so $\IS(r)\subset \Omega$, consider the map $g_r:\IS\to\IS$ defined by
 \begin{equation}
 g_r(\zeta) = \frac{f(r\zeta)}{|f(r\zeta)|}=  \frac{a \zeta^m + b \zetabar^m  + o(r^m)r^{-m}}{| a \zeta^m + b \zetabar^m  + o(r^m)r^{-m}|}, \hskip20pt |\zeta|=1.
 \end{equation}
 and 
  \begin{equation}
 g_0(\zeta) =    \frac{a \zeta^m + b \zetabar^m}{| a \zeta^m + b \zetabar^m|}, \hskip20pt |\zeta|=1.
 \end{equation}
 This is certainly continuous in $r$ for $r>0$,  but we need to analyse $r\to 0$.  As $f$ is a homeomorphism,  $|f(re^{i\theta})|\neq 0$.  Thus the denominator $| a \zeta^m + b \zetabar^m  + o(r^m)r^{-m}|$ is nonzero and converges uniformly to $| a \zeta^m + b \zetabar^m|\geq |a|-|b|$ which is nonzero by (\ref{13}). Hence  the convergence $g_r\to g_0$ is uniform.  This establishes a homotopy between $g_r$ and the map $g_0$.  Of course $g_r$ is homotopic (relative to  $\IC\setminus\{0\}$) to $f|\IS(r)\to\IC\setminus\{0\}$.
 The map $g_0$ has degree $m$.  This is a contradiction unless $m=1$ and $f_z(0)\neq 0$,   since $f|\IS(r)$ and hence $g_r$ has degree $1$. \hfill $\Box$

\section{Smoothness of $|\mu_f|^2$}
Here we establish the following theorem as a preliminary to establishing $\mu\in C^{\infty}(\Omega,\IC)$.
\begin{theorem}\label{musmooth}  Suppose that $f:\Omega \to (\tOmega,\rho)$ is harmonic.  Then $|\mu_f|^2 \in C^2(\Omega,\IR)$.
\end{theorem}
\noindent{\bf Proof.} We may assume by pre and post-composing $f$ with translations that $f(0)=0$ and the problem reduces to showing $|\mu_f|^2 \in C^2(V)$ for some neighbourhood of $0$.   If $f_z(0)\neq 0$,  then $\mu_f(z)=f_\zbar(z)/f_z(z)\in C^{\infty}(V)$,  where $V=\Omega \setminus\{z_0:f_z(z_0)=0\}$.  It certainly follows that $|\mu_f|^2$ is smooth near $0$.  Thus we need to consider the case where $f_z(0)=0$.

Thus suppose that $f_z(0)=0$.  As $J(z,f)\geq 0$,  $f_\zbar(0)=0$.   The tension equation gives us
\begin{equation}
f_{z\zbar}+a(z) f_z f_\zbar=0,\hskip20pt a(z) = (\log \rho)_w(f(z)) 
\end{equation}
and hence $f_{z\zbar}(0)=0$. 

\subsection{A first example} Let us first consider the case $f_{zz}(0)\neq 0$ (corresponding to $N=2$ below).  We have $f_z(0)=f_\zbar(0)=f_{z\zbar}(0)=0$ and $a= f_{zz}(0)\neq0$ and $b= f_{\zbar\zbar}(0)$ with $|b|\leq |a|$,  since $J(z,f)\geq 0$.  Differentiating the tension equation shows $f_{zz\zbar}(0)=f_{z\zbar\zbar}(0)=0$ and thus
\begin{eqnarray*}
f(z) & = & a  z^2 + b \zbar^2 + c z^{3} + d  \zbar^{3 }  
 + \sum_{i+j=4} a_{i,j} z^i\zbar^j +   \sum_{i+j=5} b_{i,j} z^i\zbar^j + E(z)  
\end{eqnarray*}
with $E(z)$ smooth and $|E(z)|=O(|z|^6)$.  Then
\begin{eqnarray*}
f_z(z) & = &2 a  z  + 3  c z^{2} + \sum_{i+j=4} a_{i,j} i z^{i-1} \zbar^j  + E_1(z)  \\
f_\zbar(z) & = &2 b \zbar    + 3  d \zbar^{2} + \sum_{i+j=4} a_{i,j} j z^{i} \zbar^{j-1}  + E_2(z)  
\end{eqnarray*}
with $E_k(z)$ smooth,  $k=1,2$ and $|E_k(z)|=O(|z|^4)$.
Next we see
  \begin{eqnarray*}
f_z(z) f_\zbar (z) & = &4 ab  z\zbar  + 6  bc \zbar z^{2} + 6  ad \zbar^2 z + \mbox{terms of degree $\geq 4$}  \\
f_{z\zbar}(z) & = &  \sum_{i+j=4 } a_{i,j} i j z^{i-1} \zbar^{j-1} + \sum_{i+j=5 } b_{i,j} ij z^{i-1}\zbar^{j-1} + \mbox{terms of degree $\geq 4$}  
\end{eqnarray*}
  which gives $a_{i,j}=0$ unless  $i=j=2$,  $i=4$ or $j=4$.  We therefore obtain the expansion
  \begin{eqnarray*}
f(z) & = & a  z^2 + b \zbar^2 + c z^{3} + d  \zbar^{3 }  + e z^4 + d \zbar^4 + g z^2 \zbar^2 + E(z)  
\end{eqnarray*}
  and $E(z)=O(|z|^5)$.  From this we calculate that
    \begin{eqnarray*}
f_z(z) & = & 2a  z+3c z^{2} +4 e z^3 + 2 g z  \zbar^2 + E_1(z)  \\
& = & z(2a   +3c z  +4 e z^2 + 2 g   \zbar^2 + \tilde{E}_1(z) ) = z u(z)\\
f_{\zbar}(z) & = & 2 b \zbar  + 3 d  \zbar^{2 }  +4 d \zbar^3 + 2 g z^2 \zbar  + E_2(z) \\
& = & \zbar( 2 b + 3 d  \zbar   +4 d \zbar^2 + 2 g z^2  + \tilde{E}_2(z) ) = \zbar v(z)
\end{eqnarray*}
and $\tilde{E}_i=O(|z|^3)$,  $i=1,2$.
Now $u$ and $v$ are $C^2$ functions near $0$,  but they might not have any higher degree of differentiability.  For instance the error term in $u$ might contain terms such as $\zbar^4/z$,  so $u_{\zbar\zbar\zbar}$ may not be continuous.  In any case,  we now have
\begin{equation}
|\mu_f(z)|^2 =\frac{|f_\zbar(z)|^2}{|f_z(z)|^2} =  \frac{|v(z)|^2}{|u(z)|^2}
\end{equation}
with $u,v \in C^2$ near $0$ and $u\neq 0$.  Hence $|\mu_f|^2\in C^2$ near $0$.
\subsection{Higher order analysis}
We now have to extend the above argument should $f_{zz}(0)=0$.  Indeed,  suppose
\begin{equation}\label{9} \Big(\frac{\partial^m\; f}{\partial z^m}\Big)(0) =  \Big(\frac{\partial^m\; f}{\partial \zbar^m}\Big)(0) = 0
\end{equation}
 for all $m\leq N$,  then an easy induction,  repeatedly differentiating the tension equation,  implies that for all $p$ and $q$,  $p+q=N+1$
  \[ \Big(\frac{\partial^m\; f}{\partial z^p \partial \zbar^q}\Big)(0) = 0.\]
 Since $f$ can only have a zero of finite order,  for otherwise Hopf differential is identically $0$ and $f$ is conformal, there is an $N$ such that
 \begin{equation}
a_{N,0} =  \Big(\frac{\partial^m\; f}{\partial z^m}\Big)(0) \neq 0,  \hskip10pt    a_{0,N}=\Big(\frac{\partial^m\; f}{\partial \zbar^m}\Big)(0), \hskip5pt{\rm and} \hskip10pt |a_{N,0}|\geq |a_{0,N}|
 \end{equation}
The last inequality follows since  the smooth function $J(z,f)$ is non-negative. Further
  \begin{equation}
 \Big(\frac{\partial^m\; f}{\partial z^p \partial \zbar^q}\Big)(0) = 0, \hskip5pt \mbox{whenever $i+j\leq N+1$ and $i,j\geq 1$.}
  \end{equation}
This gives us the following form for the Taylor series of $f$
\begin{eqnarray*}
f(z) & = & a_{N,0} z^N + a_{0,N} \zbar^N + a_{N+1,0} z^{N+1} + a_{0,N+1} \zbar^{N+1}  \\
& & + \sum_{i+j=N+2} a_{i,j} z^i\zbar^j + \sum_{i+j=N+3} b_{i,j} z^i\zbar^j +E(z)  
\end{eqnarray*}
  where the error term $E(z)=O(|z|^{N+3})$ and $N\geq 3$.    We follow the argument above to compute
  \begin{eqnarray*}
f_z(z) & = & \alpha z^{N-1} +  \beta  z^{N}  + \sum_{i+j=N+2} i a_{i,j} z^{i-1} \zbar^j + \sum_{i+j=N+3} i b_{i,j} z^{i-1} \zbar^j +E_1(z)  \\
f_\zbar(z) & = & \gamma \zbar^{N-1} + \delta \zbar^{N}   + \sum_{i+j=N+2} a_{i,j} j z^i\zbar^{j-1} + \sum_{i+j=N+3}j  b_{i,j}  z^i\zbar^{j-1} +E_2(z)  
\end{eqnarray*}
The tension equation then implies that the lowest degree term of $f_{z\zbar}$ is of order $2N-2$ and in fact 
  \begin{eqnarray*}
  f_{z\zbar}(z) & = &  a  z^{N-1}\zbar^{N-1} + b  z^{N-1}\zbar^{N} + c z^{N}\zbar^{N-1} + \mbox{terms of degree 2N or more}
\end{eqnarray*}
Possibly some of these coefficients are zero of course,  depending on the expansion of $(\log\rho)_w(f)$.
Thus $f$ has the form
\begin{eqnarray*}
f(z) & = & z^N P(z) + \zbar^N Q(\zbar) +    a z^N\zbar^N +  b  z^{N }\zbar^{N+1} + c z^{N+1}\zbar^{N} \\ && + \mbox{terms of degree 2N+2 or more}
\end{eqnarray*}
Where $P$ and $Q$ are polynomials of degree $N+1$ with $P(0)\neq0$.  Hence
\begin{eqnarray*}
f_z(z) & = & N z^{N-1} P(z)+z^NP'(z) +  aN z^{N-1}\zbar^N +  b  z^{N-1 }\zbar^{N+1} + c (N+1) z^{N}\zbar^{N}  \\ && + \mbox{terms of degree 2N+1 or more}\\
& = & z^{N-1} U(z)  
\end{eqnarray*}
where $U(z)=N   P(z)+z P'(z) +  aN  \zbar^N +  b \zbar^{N+1} + c (N+1) z \zbar^{N} + E(z)$ and $|E(z)|= O(|z|^{N+2})$.  Hence $U\in C^{N}$ near $0$.

A similar expansion holds for $f_\zbar$ and as with the earlier example,  the result follows. \hfill $\Box$

 \section{Sub/super-harmonicity of $\log |\mu_f|$}
 
Let $\Phi_f$ be the holomorphic Hopf differential of the harmonic map $f:\Omega \to (\tOmega,\rho)$.  Fix a point $z_0 \in \Omega\setminus \{z:f_\zbar(z) = 0 \}$.  We choose a well defined branch of the argument in a   domain containing  $f_z(z_0)$ and $f_\zbar(z_0)$ so that we can work with well defined functions $\log f_z = \log |f_z|+i \arg(f_z)$ and also $\log f_\zbar= \log |f_\zbar|+i \arg(f_\zbar)$.  Notice that as $\phi_f=\rho^2(f) f_z \overline{f_\zbar}$ is holomorphic,  $\arg(f_z \overline{f_\zbar})$ is harmonic.  This is basically why the imaginary terms disappear below.   Then,  as $f_{z\zbar}=-(\log \rho)_w(f) f_z f_\zbar$, 
   \begin{eqnarray*}
\lefteqn{ \frac{1}{4}\Big(       \Delta   \log   {  f_\zbar} - \Delta   \log   {f_z } \Big) =    \Big[\frac{ {  f_{z\zbar}}}{ {  f_{\zbar}} } \Big]_\zbar - \Big[ \frac{{f_{z\zbar} }}{{f_{z} }} \Big]_z }\\ 
  & = & \Big[  {{(\log \rho)_w(f)   f_\zbar }} \Big]_z - \Big[  { (\log \rho)_w(f) f_z  } \Big]_\zbar\\
   & = &  
  \Big[  {(\log \rho)_w(f)\Big]_z   {f_\zbar }} + {(\log \rho)_w(f)  {f_{z\zbar} }}  -  \Big[   (\log \rho)_w (f) \Big]_\zbar {  f_z  }  -(\log \rho)_w (f)   {  f_{z\zbar}  }  \\
& = &    \Big[  {(\log \rho)_w(f)\Big]_z    {f_\zbar }}  -  \Big[   (\log \rho)_w (f) \Big]_\zbar {  f_z  }    \\
& = &    \Big[  {(\log \rho)_{ww}(f) f_z + (\log \rho)_{w\wbar}(f)\overline{f_\zbar} \Big]    {f_\zbar }}  -  \Big[   (\log \rho)_{ww} (f) f_\zbar + (\log \rho)_{w\wbar} (f) \overline{f_z} \Big]{  f_z  }    \\
& = &    (\log \rho)_{w\wbar} (f)| f_\zbar |^2     - (\log \rho)_{z\zbar} (f)|f_z|^2 =   -   \frac{1}{4} (\Delta \log \rho)(f) J(z,f).
 \end{eqnarray*}
 We can now prove the following theorem.
 \begin{theorem}\label{thm5} Let $f:\Omega\to (\tOmega,\rho)$ be harmonic.  Then,  away from $\{\mu_f = 0 \}$, 
 \begin{equation}\label{musub}
 \Delta \log |\mu_f| = -  (\Delta \log \rho)(f) J(z,f)
 \end{equation}
 \end{theorem}
 \noindent{\bf Proof.}  We evidently have the equality (\ref{musub}) away from the discrete set where $f_\zbar =0$ from the above calculation.  However,  by Theorem \ref{musmooth}, $|\mu_f|^2$,  and hence $\log|\mu_f|$,  is twice continuously differentiable and so the left-hand side of (\ref{musub}) is well defined and continuous  on $\Omega\setminus \{\mu = 0\}$.  The right hand side is continuous in $\Omega$ and so the result follows.  \hfill $\Box$

 \bigskip
 
The curvature of the metric $\rho(z)|dz|$ is 
 $$ {\cal K}_\rho (w) = \frac{-(\Delta \log \rho)(w)}{\rho^2(w)} $$
This gives us the following version of Theorem \ref{thm5}

\begin{theorem}\label{thm6} Let $f:\Omega\to (\tOmega,\rho)$ be harmonic.  Then,  away from $\{\mu = 0 \}$, 
 \begin{equation}\label{22}
 \Delta \log |\mu_f| = {\cal K}_\rho(f) \rho^2(f) J(z,f)
 \end{equation}
 \end{theorem}
 This immediately shows $\log |\mu_f| $ is subharmonic for a non-negatively curved metric and superharmonic for a non-positively curved metric.  Since subharmonic functions cannot achieve a local maximum at an interior point,  we immediately see $\log |\mu_f| < 0$,  that is locally $|\mu_f|<1$ and $f$ is locally quasiregular.   

  \subsection{Super averaging property.}
 
In order to establish more generally that $|\mu_f|<1$ locally,  we need an analysis near the points where potentially $|\mu_f|=1$.   We assume that $0\in \Omega$ is such a point,  $|\mu_f(0)|$ is close to $1$ (possibly equal).  In a neighbourhood of $0$ we have
 \[ w(z) \log \frac{1}{|\mu_f|} = (1-|\mu_f|^2) \]
 for $w(z)$ a smooth positive real valued function,  $w(0)=1$.  Then,  near $0$, 
\begin{eqnarray*} 
{\cal K}_\rho(f) \rho^2(f) J(z,f) & = & {\cal K}_\rho(f) \rho^2(f) (1-|\mu_f|^2) |f_z|^2 \\
& = & w(z) {\cal K}_\rho(f) \rho^2(f) |f_z|^2 \;  \log  \frac{1}{|\mu_f|}  
\end{eqnarray*}
 Next,  with $u(z)=-\log |\mu_f(z)| \geq 0$,  (\ref{22}) reads as
 \begin{equation}
  \Delta u = \big( - a(z) {\cal K}_\rho(f) \rho^2(f) |f_z|^2\big) \, u
 \end{equation}
 The term in brackets is  continuous,  so bounded near $0$,  and we obtain for some $C\geq 0$, 
  \begin{equation}
  \Delta u \leq C u
 \end{equation}
 It seems E. Heinz \cite{Heinz} was the first to point out a super-averaging  principle for these subsolutions to the eigenvalue equation of the Laplacian.  Rather more sophisticated versions of this are given in a  paper of Choi and Treibergs, \cite{CT}.  Here we can give a sharpened version of these results in our setting.
 
 \begin{lemma}\label{Heinz}  Let $u$ be a $C^2(U)$ non-negative real valued function defined on a neighbourhood $U$ of $0\in \IC$.  Let $d={\rm dist}(0,\partial U)$ and suppose that there is a constant $C \geq 0$ so that $\Delta u \leq C \, u$ on $ U$.   Then 
 \begin{equation}\label{lem1}
 u(0) \geq \frac{1}{2\alpha^2} \;  \int_{\ID(\alpha)} u(z) \; dz
 \end{equation}
where  $\alpha =  \min \{\frac{1}{4} \sqrt{\frac{e}{C}},\frac{d}{2} \}>0$.
 Consequently if $u(z)>0$ almost everywhere,  we have $u(0) > 0$.
 \end{lemma}
 \noindent{\bf Proof.}  Let $v=u(d\, z)$ so $v$ is defined on the unit disk $\ID$ and 
 \[ (\Delta v)(z) = d^2( \Delta u)(d\, z) \leq C d^2  u(d\,z) = Cd^2 v(z).\]
  We multiply by $\log |z|$,  the fundamental solution to the Laplacian - meaning that as distributions $\Delta \log |z| = \delta_0$, the Dirac delta function.   We then integrate by parts twice (omitting the details of the standard distributional calculation to achieve (\ref{diet})),  then with the   outer-normal $\eta=z/|z|$ to the circle we see
 \begin{eqnarray}
 \int_\ID \Delta v  \, \log|z| \; dz  & \geq & Cd^2 \int_\ID v \log|z| \; dz \nonumber \\
  \int_\ID   v  \, \Delta \log|z| \; dz  & -   &  \oint_{|z|=1} v (\nabla \log |z| )\cdot\eta + \oint_{|z|=1}  \log|z| \nabla v \cdot\eta   \label{diet} \\  & \geq &  Cd^2 \int_\ID v \log|z| \; dz \nonumber
\end{eqnarray}
Note that $\int_\ID v \log|z| \; dz\leq 0$.  Hence, as $v(0)=u(0)$, a little calculation yields
\begin{equation}\label{25}
u(0) \geq    \oint_{|z|=1} v      +    Cd^2 \int_\ID v \log|z| \; dz  
 \end{equation}
 Now,    we can again follow the above argument with $v$ replaced by $v_s(z) = v(s z)$ to see,  as $\Delta v_s \leq C d^2 s^2 \, v_s $,   that we obtain (\ref{25})  in the form
  \begin{equation}\label{26}
u(0) \geq    \oint_{|z|=1} v_s(z)      +    Cd^2s^2 \int_\ID v_s(z) \log|z| \; dz 
 \end{equation}
 Then multiplying by $s$ and integrating from $0$ to $\delta<1$ yields
  \begin{eqnarray*}
\int_{0}^{\delta} u(0)\,  s \, ds & \geq &  \int_{0}^{\delta}  \oint_{|z|=1} v(sz) \, s\, ds      +    Cd^2  \int_\ID \int_{0}^{\delta} v(s z) s^3 \, ds \log|z| \; dz \\
& \geq &  \int_{\ID(\delta)}   v(z) \, dz      +    Cd^2  \int_{0}^{1} \int_{0}^{\delta} \oint_{|z|=1} v(r s z) s^3 \, ds d \theta \, r \log r  \, dr \\
& \geq &  \int_{\ID(\delta)}   v(z) \, dz      +    Cd^2  \int_{0}^{1} \int_{\ID(\delta)} v(r z) |z|^2 dz  \, r \log r  \, dr \\
& = &  \int_{\ID(\delta)}   v(z) \, dz      +    Cd^2  \int_{0}^{1} \int_{\ID(r\delta)} v(w) |w|^2 dw \frac{\log r }{r} \, dr \\
& = &  \int_{0}^{1} \Big(\int_{\ID(\delta)}   v(z) \, dz      +    Cd^2   \int_{\ID(r\delta)} v(w) |w|^2 dw \frac{\log r }{r} \Big) dr \\
& \geq &  \int_{0}^{1} \Big(\int_{\ID(r\delta)}   v(w) \, dw      +    Cd^2   \int_{\ID(r\delta)} v(w) |w|^2 dw \frac{\log r }{r} \Big) dr \\
& \geq &  \int_{0}^{1}  \int_{\ID(r\delta)}   v(w) \Big( 1      +    Cd^2  |w|^2  \frac{\log r }{r} \Big) dw\, dr \\
& \geq &  \int_{0}^{1}  \int_{\ID(r\delta)}   v(w) \Big( 1      +    \delta^2 Cd^2  r \,  \log r   \Big) dw\, dr 
 \end{eqnarray*}
 Since $r\log r \geq -1/e$,  if we choose 
$\delta  = \min \{\frac{1}{2d} \sqrt{\frac{e}{C}},1\}$,
 then we achieve
   \begin{eqnarray*}
u(0) & \geq & \frac{1}{\delta^2} \int_{0}^{1}  \int_{\ID(r\delta)}   v(w) dw\, dr 
 \end{eqnarray*}
Now $v(w)=u(d\; w)\geq 0$ and the monotonicty of the integrand $ \int_{\ID(r\delta)}   v(w) dw$ gives
    \begin{eqnarray*}
u(0) & \geq & \frac{1}{\delta^2} \int_{1/2}^{1}  \int_{\ID( \delta/2)}  u(dw)  \, dr    =  \frac{1}{2d^2\delta^2}   \int_{\ID( d \delta/2)}  u(z)  dz  
 \end{eqnarray*}
 and the result follows. \hfill $\Box$
 \bigskip
 
 \subsection{Completion of the proof.}
 We can now apply Lemma \ref{Heinz} to $u=-\log|\mu_f|$ which is certainly positive almost everywhere,  to see that at each point of $\Omega$ we have $|\mu_f|<1$ and $|\mu_f|\in C^2(\Omega)$.  Thus,   on any relatively compact $U\subset\Omega$,  we have $|\mu_f|\leq k_U < 1$ and $f$ is therefore quasiregular.  The Stoilow factorisation theorem \cite[Corollary14.4.5]{AIM} asserts that there is a quasiconformal $g:U\to \IC$ and a holomorphic $\varphi:g(U)\to \IC$ with 
$f|U = \varphi \circ g$.  If $f$ has local degree $1$,  then so does $\varphi$ and $f|U$ is a local diffeomorphism by Theorem \ref{hqc}.  This establishes the main result,  Theorem \ref{MR}.

\bigskip

\noindent G.J. Martin \\ Massey University,  Auckland, NZ and  \\ Magdelen College,  Oxford. \\g.j.martin@massey.ac.nz

\end{document}